\documentclass[11pt]{amsproc}
\usepackage{amsmath,amsthm,anysize,amssymb,centernot, nicematrix, framed, tikz,comment}
\usepackage{booktabs, float, lmodern}
\usepackage[dvipsnames,x11names]{xcolor}
\usepackage{hyperref}
\hypersetup{
    colorlinks,
    linkcolor={red!60!black},
    citecolor={green!60!black},
    urlcolor={blue!60!black}
}

\usepackage{enumerate}
\usepackage{url}
\usepackage[T1]{fontenc}
\usepackage{geometry}
\theoremstyle{plain}
\newtheorem{Theorem} {Theorem} [section]
\newtheorem{Proposition} [Theorem] {Proposition}
\newtheorem{Lemma} [Theorem] {Lemma}
\newtheorem{Corollary} [Theorem] {Corollary}
\newtheorem{Example} [Theorem] {Example}

\newtheorem{Conjecture} [Theorem] {Conjecture}

\newtheorem{Question} [Theorem] {Question}

\theoremstyle{definition}
\newtheorem{Definition} [Theorem] {Definition}

\DeclareMathOperator{\F}{\mathbb{F}}

\DeclareMathOperator{\GL}{\mathrm{GL}}

\DeclareMathOperator{\PSL}{\mathrm{PSL}}

\DeclareMathOperator{\diag}{\mathrm{diag}}

\DeclareMathOperator{\abs}{\mathrm{abs}}

\renewcommand{\mp}{\mathcal P}
\newcommand{\ml}{\mathcal L}
\newcommand{\mf}{\mathcal F}

\DeclareMathOperator{\MOAS}{\mathrm{MOAS}}

\DeclareMathOperator{\rank}{\mathrm{rank}}

\renewcommand{\le}{\leqslant}
\renewcommand{\ge}{\geqslant}

\title{On configurations that are nearly generalised polygons}
\author{Sam Adriaensen}
\address[Adriaensen]{Department of Mathematics and Data Science, 
Vrije Universiteit Brussel, 
Pleinlaan 2, 1050 Elsene, Belgium.}
\email{sam.adriaensen@vub.be}

\author{John Bamberg}
\address[Bamberg]{Department of Mathematics and Statistics, 
The University of Western Australia, 
35 Stirling Highway, 
Perth, W.A. 6009, Australia.}
\email{John.Bamberg@uwa.edu.au}

\author{Jason Williford}
\address[Williford]{Department of Mathematics and Statistics,
University of Wyoming,
1000 E. University Ave.,
Laramie, WY 82071 USA}
\email{jwillif1@uwyo.edu}



\begin{document}

\begin{abstract}
Payne (1972) introduced nearly generalised $n$-gons as $(s+1)$-regular bipartite graphs of girth $2n$ with small excess over the Moore bound. 
Biggs and Ito (1980) proved nonexistence of such graphs if $s > 1$ and $n > 3$, leaving open the existence of regular antipodal planes.
We revisit this problem from several equivalent viewpoints.  First we show that a regular antipodal plane of order $s$ is equivalent to a balanced weighing matrix whose support is the incidence matrix of a biplane with parameters
\[
        2-\left(\frac{s^2+s+2}{2},\,s+1,\,2\right).
\]
This gives short determinant and local-arithmetic obstructions to existence.  We then formulate the question of signing a biplane as an affine linear system over $\mathbb F_2$, and use this formulation to prove that no regular antipodal plane of order $s$ exists with $4\le s\le 14$.
Finally, we give further equivalent formulations in terms of sharply deranged families of permutations and mutually orthogonal array-bic squares. These translations provide new computational and structural approaches to Payne's conjecture that regular antipodal planes exist only for $s=2,3$.
\end{abstract}

\maketitle

\section{Introduction}

The \emph{cage problem} is to determine the minimum number $c(k,g)$ of vertices in a $k$-regular graph of girth $g$; graphs attaining this bound are called \emph{cages}. Although cages exist for all integers $k \ge 2$ and $g \ge 3$, the exact value of $c(k,g)$ is known only for a small number of parameter pairs, together with three infinite families occurring when $g \in \{6,8,12\}$ and $k-1$ is a prime power. These families arise as the incidence graphs of generalised polygons. For the even-girth case, by counting out from a vertex, it is not difficult to see that a simple undirected $(s+1)$-regular bipartite graph of even girth $g$ has at least
\[
2\left(1+s+s^2+\cdots+s^{(g-2)/2}\right)
\]
vertices, and those graphs meeting this bound are known as \emph{Moore graphs}. Interestingly, there are very few
Moore graphs. It has been proved that a bipartite Moore graph
with $k>2$ has $g \in \{4, 6, 8, 12\}$ and is the incidence graph of a generalised $(g/2)$-gon of order $k-1$.

Payne (1972) \cite{Payne:1972aa} introduced the notion of a \emph{nearly} generalised $n$-gon, by relaxing
the condition that the number of vertices meets the Moore bound. These graphs have parameters $(n,s)$, and they are $(s+1)$-regular bipartite graphs of girth $2n$ with
\[
2\left(\mathbf{2}+s+s^2+\cdots+s^{n-1}\right)
\]
vertices. So for example, if $n=3$, it is the incidence graph of a \emph{regular antipodal plane}:
a partial linear space of points and lines such that (i) every point is incident with $s+1$ lines,
(ii) every line is incident with $s+1$ points, and (iii) for every point $P$, there is a unique
point $P^*$ not collinear with $P$, and for every line $\ell$, there is a unique line $\ell^*$
not concurrent with $\ell$. Antipodal planes were studied by Dembowski \cite[Section 7.4]{Dembowski:1997aa} 
as elliptic semi-planes with parallel classes of size 2, and he observed that they are 2-covers
of biplanes. Payne showed that there are examples of nearly 3-gons for $s=2,3$, but no examples
for $s=4$. Indeed, he conjectured that there are no examples for $s>3$.
He also showed that there are no nearly $n$-gons for $n\in \{4,5,6\}$ (with $s>1$).
A substantial breakthrough was made by Biggs and Ito \cite{BiggsIto} in 1980: there is no nearly $n$-gon
(with $s>1$) for $n>3$. So the problem reduces to the existence of regular antipodal planes.
Throughout, we will assume $s>1$.

Further work has been done recently by Schleiermacher \cite{Schleiermacher:2016aa}
and De Boeck and Van de Voorde \cite{De-Boeck:2022aa}. The latter showed that no antipodal plane of order $s>3$ can be embedded in a Desarguesian projective plane. Both papers note that a regular antipodal plane double covers a biplane, but did not note the existence
of a weighing matrix, which imposes further restrictions on the possible orders of a regular antipodal plane.

In this paper we draw attention to some objects in graph theory, incidence geometry, and design theory that can be unified.
For instance, a \emph{balanced weighing matrix} with parameters $(k,2)$ is equivalent to a regular antipodal plane
of order $k-1$ (see Theorem \ref{thm:BWM_equivalence}). Figure \ref{fig:3equivalent} gives an example
of the three combinatorial objects that we will show are equivalent. Indeed, a regular antipodal plane is equivalent to a \emph{sharply deranged} set of permutations, and hence gives rise to an \emph{equidistant permutation array}. Through this connection, a recent result of Tolonen \cite{Tolonen} implies that there is no regular antipodal plane of order 8. We show that there are no regular antipodal planes of order $s$ with $4 \le s \le 14$ (see Section \ref{section:3gons}).

\begin{figure}[!ht]
\begin{tabular}{c|c|c}
$
\begin{aligned}
&(1\,2\,3\,4), (1\,2\,4\,3), (1\,3\,2\,4),\\
&(1\,3\,4\,2), (1\,4\,2\,3), (1\,4\,3\,2)
\end{aligned}$
&
$
\begin{bmatrix}
\cdot&3&2\\
0&\cdot&3\\
0&2&\cdot\\
\end{bmatrix}\,
\begin{bmatrix}
\cdot&0&3\\
3&\cdot&1\\
1&0&\cdot\\
\end{bmatrix}\,
\begin{bmatrix}
\cdot&2&0\\
1&\cdot&0\\
2&1&\cdot\\
\end{bmatrix}
$
&
${\tiny \begin{bmatrix}
0& +& 0& +& +& 0& +\\
0& +& +& 0& 0& +&-\\
0& 0& +& +&-&-& 0\\
+& 0& 0& 0& +&-&-\\
+& 0& +&-& 0& 0& +\\
+&-& 0& +& 0& +& 0\\
+& +&-& 0&-& 0& 0
\end{bmatrix}}$\\[2em]
Sharply deranged set & Mutually orth.\ array-bic squares & Balanced weighing matrix
\end{tabular}
\caption{Three equivalent objects.}\label{fig:3equivalent}
\end{figure}

\subsection*{Declaration of AI use}
We used AI (ChatGPT) to find references, such as the important work of Biggs and Ito \cite{BiggsIto}
and the result of Schellenberg \cite{Schellenberg1975}. 
It also gave us the main idea of Proposition \ref{exact-rank}. We 
used it to write code to search for array-bic squares in \textsf{Julia} \cite{Julia}.

\section{Nearly 3-gons, antipodal planes}\label{section:3gons}

\subsection{Weighing matrices}

A regular antipodal plane has the property that the star operation preserves incidence:
if $P$ is a point incident with $\ell$, then $P^*$ is incident with $\ell^*$ (see \cite[Lemma 2.2]{De-Boeck:2022aa}).
We can define a biplane, a symmetric $2-((s^2+s+2)/2,s+1,2)$ design, in the following way:
the points are the pairs $\{P,P^*\}$, and the blocks are the pairs $\{\ell,\ell^*\}$,
and a point $\{P,P^*\}$ is incident with a block $\{\ell,\ell^*\}$ if 
$P$ is incident with $\ell$ or $P^*$ is incident with $\ell$.
Payne observes that the Bruck-Ryser-Chowla Theorem for symmetric designs,
plus a linear algebraic property of the incidence matrix, give rise to infinitely many
forbidden values of $s$:

\begin{Theorem}[Payne {\cite[Theorems 5.1, 5.2]{Payne:1972aa}}]\leavevmode
\begin{enumerate}[(a)]
    \item 
If $s\equiv 3\pmod{8}$ and if there is a prime congruent to $7\pmod{8}$
dividing the square-free part of $s-1$, then no nearly 3-gon of order $s$ exists.
    \item 
If $s\equiv 1\pmod{8}$ and if there is a prime congruent to $3\pmod{8}$
dividing the square-free part of $s+1$, then no nearly 3-gon of order $s$ exists.
\end{enumerate}
\end{Theorem}

We observe that part of this result can be proved via the theory of \emph{weighing matrices},
thus giving an alternative and shorter proof of \cite[Theorem 5.2]{Payne:1972aa}.
A \emph{weighing matrix} of order $t$ and weight $w$ is a matrix $W$ with entries in $\{-1,0,1\}$ such that $WW^\top=wI_t$. 
The absolute matrix $\abs A$ of a $\{-1,0,1\}$-matrix $A= (a_{i,j})$ is the matrix
$(a_{i,j}^2)$. A weighing matrix $W$ is a \emph{balanced weighing matrix} with parameters $(k,\lambda)$
if 
\[
(\abs W)(\abs W)^\top=(k-\lambda)I+\lambda J.
\]
Indeed, this is equivalent to $\abs W$ being the incidence matrix of a symmetric balanced incomplete block
design.

\begin{Theorem}\label{thm:BWM_equivalence}
    A regular antipodal plane of order $s$ is equivalent to a balanced weighing matrix with $k=s+1$ and $\lambda=2$.
\end{Theorem}

\begin{proof}
    First suppose that $\mathcal A$ is a regular antipodal plane of order $s$.
Let $P\mapsto P^*$ be the antipodal involution on points, and
$\ell\mapsto \ell^*$ the antipodal involution on lines. Thus
$|\mathcal P|=|\mathcal L|=s^2+s+2$,
and every point-pair $\{P,P^*\}$ and every line-pair
$\{\ell,\ell^*\}$ has size $2$. Hence there are
$v=\frac{s^2+s+2}{2}$
point-pairs and the same number of line-pairs.
Choose one representative $P_i$ from each point-pair and one representative
$\ell_a$ from each line-pair. Define a $v\times v$ matrix $W$, with rows
indexed by line-pairs and columns indexed by point-pairs, by
\[
        W_{ai}=
        \begin{cases}
        1, &\text{if } P_i\in \ell_a,\\
        -1,&\text{if } P_i^*\in \ell_a,\\
        0, &\text{if neither } P_i \text{ nor } P_i^* \text{ lies on } \ell_a.
        \end{cases}
\]
A line contains exactly $s+1$ points, no two of which are antipodal, so each
row of $W$ has exactly $s+1$ nonzero entries.

We claim that distinct rows of $W$ are orthogonal. Let $\ell$ and $m$ be
representatives of two distinct line-pairs. By regularity,
$\ell$ meets both $m$ and $m^*$ in exactly one point. Write
$\ell\cap m=\{P\}$ and $\ell\cap m^*=\{Q\}$. Then the corresponding two point-pairs $\{P,P^*\}$ and $\{Q,Q^*\}$ are
distinct. In the column corresponding to $\{P,P^*\}$, the two rows have the
same sign, since both $\ell$ and $m$ contain the same point $P$. Hence
this column contributes $+1$ to the row inner product. In the column
corresponding to $\{Q,Q^*\}$, the row for $\ell$ contains $Q$, whereas
the row for $m$ contains $Q^*$, because $m^*$ contains $Q$. Hence this
column contributes $-1$. All other columns contribute $0$. Therefore the
inner product of the two rows is $1-1=0$.
Thus
\[
        WW^{\top}=(s+1)I_v.
\]

Now let $B:=\abs W$. Each row of $B$ has size $s+1$.
Moreover, two distinct line-pairs meet in exactly two point-pairs, namely the
point-pairs containing $\ell\cap m$ and $\ell\cap m^*$. Hence
\[
        BB^{\top}
        =
        (s+1)I_v+2(J_v-I_v)
        =
        (s-1)I_v+2J_v.
\]
Equivalently, $B$ is the incidence matrix of a symmetric
$2$-$(v,s+1,2)$ design; a biplane. Therefore $W$ is a balanced weighing matrix with
$k=s+1$ and $\lambda=2$.

Conversely, suppose that $W$ is a balanced weighing matrix with parameters
$(s+1,2)$. Thus $WW^{\top}=(s+1)I_v$, and $B=\abs W$ is the incidence matrix of a biplane. 
The equation $\lambda(v-1)=k(k-1)$ for symmetric designs gives
$2(v-1)=(s+1)s$, so $2v=s^2+s+2$. We construct an antipodal plane. For each column $i$ of $W$, introduce two
points $P_i^+$ and $P_i^-$, and declare them antipodal. For each row $a$, introduce two lines
$\ell_a^+$ and $\ell_a^-$, and declare them antipodal. Define incidence by
\[
        P_i^\epsilon \in \ell_a^\delta
        \quad\Longleftrightarrow\quad
        W_{ai}=\delta\epsilon,
        \qquad \epsilon,\delta\in\{\pm1\}.
\]
Thus if $W_{ai}=0$, neither $P_i^+$ nor $P_i^-$ lies on either
$\ell_a^+$ or $\ell_a^-$; while if $W_{ai}\neq0$, each of
$\ell_a^+$ and $\ell_a^-$ contains exactly one point from the antipodal
pair $\{P_i^+,P_i^-\}$.
Since each row of $W$ has exactly $s+1$ nonzero entries, each line contains
$s+1$ points. Similarly, since $W$ is square and $WW^{\top}=(s+1)I_v$,
we also have $W^{\top}W=(s+1)I_v$, so each column has $s+1$ nonzero
entries; hence each point lies on $s+1$ lines. The total number of points and
lines is $2v=s^2+s+2$.
It remains to check the partial linear space condition and the antipodal property.
The two points $P_i^+$ and $P_i^-$ are not collinear, because every line
contains at most one of them by construction.
Now take two non-antipodal points $P_i^\epsilon$ and $P_j^\eta$, with
$i\neq j$. Since $B$ is a symmetric $2$-$(v,s+1,2)$ design, the two
columns $i$ and $j$ have common support in exactly two rows, say $a$ and
$b$. Since $W^{\top}W=(s+1)I_v$, distinct columns of $W$ are orthogonal,
so
\[
        W_{ai}W_{aj}+W_{bi}W_{bj}=0.
\]
Thus one of these two products is $+1$, and the other is $-1$. Therefore
exactly one of the two rows, say $a$, satisfies
\[
        W_{ai}W_{aj}=\epsilon\eta.
\]
For this row there is a unique choice of sign $\delta\in\{\pm1\}$ such that
\[
        W_{ai}=\delta\epsilon,
        \qquad
        W_{aj}=\delta\eta.
\]
Hence $P_i^\epsilon$ and $P_j^\eta$ lie on the unique line
$\ell_a^\delta$. Therefore any two non-antipodal points are joined by exactly
one line, while antipodal points are joined by no line.
This is precisely a regular antipodal plane of order $s$.

Finally, the two constructions are inverse up to the choices of
representatives. Changing the representative of a point-pair multiplies the
corresponding column of $W$ by $-1$, and changing the representative of a
line-pair multiplies the corresponding row by $-1$. Thus the equivalence is
between regular antipodal planes of order $s$ and balanced weighing matrices
with parameters $k=s+1$, $\lambda=2$,
up to the usual row and column sign changes and permutations.
\end{proof}

\begin{Corollary}\label{corollary:newcondition}
    If a regular antipodal plane exists of order $s\equiv 0,3\pmod{4}$, then $s+1$ is a square.
\end{Corollary}

\begin{proof}
Let $v:=\tfrac{1}{2}(s^2+s+2)$, the order of $W$. Now
    \[
    \det(WW^\top)=\det(W)^2=(s+1)^v.
    \]
So $(s+1)^v$ is a square. Hence, if $v$ is odd, then $s+1$ is a square. Note that $v$ is odd precisely when $s\equiv 0,3\pmod{4}$.
\end{proof}

Before we show a unified approach to showing that certain families of weighing matrices do not exist, we address an
observation by Mullin \cite{Mullin:1978aa}. He lists necessary conditions for the existence of a balanced weighing matrix with parameters $(k,\lambda)$. The last one is attributed to Raghavarao \cite{Raghavarao}, which we reproduce below. 
However, Raghavarao (in Section 4) does not explicitly state Mullin's general condition (iv). He proves the 
`$W(N,N-1)$ case', and his proof immediately extends to general $k$. 

Let $X:=BW^\top$ and $r=s+1$.
Since $\det(XX^\top)=r^{v+2}(r-2)^{v-1}$ is a square, we recover the biplane and weighing matrix conditions on $r$ and $r-2$ in one elegant form.
Set $G:=XX^{\top}=r(r-2)I+2rJ$.
Since $X$ is an integral nonsingular matrix, $\det(G)=(\det X)^2$ is a square. The eigenvalues of $G$ are
$r^3$ with multiplicity $1$, and $r(r-2)$ with multiplicity $v-1$. Hence
\[
\det(G)=r^3\{r(r-2)\}^{v-1}
       =r^{v+2}(r-2)^{v-1}.
\]
Therefore $r^{v+2}(r-2)^{v-1}$ is a square. It follows that
\[
\begin{cases}
r \text{ is a square}, & v \text{ odd},\\[2mm]
r-2 \text{ is a square}, & v \text{ even}.
\end{cases}
\]
In particular, $r$ cannot be congruent to $5$ or $7$ modulo $8$.

We now derive Raghavarao's condition using the Hasse-Minkowski theorem.
Let $p$ be an odd prime, and $x$ and $y$ rational numbers.
The \emph{Hilbert symbol} $(x,y)_p$ is defined as follows (see \cite[\S III.1.2, Theorem 1]{Serre}). Let $\alpha$ and $\beta$ be the $p$-adic valuations of $x$ and $y$ respectively, and write $x = p^\alpha u$ and $y = p^\beta v$. Then
\[
 (x,y)_p = (-1)^{\alpha \beta \frac{p-1}2} \left( \frac u p \right)^\beta \left( \frac v p \right)^\alpha, 
\]
where $\left( \frac u p \right)$ and $\left( \frac v p \right)$ are Legendre symbols.
From \cite[\S IV.3.3]{Serre}, one can derive the following version of the Hasse-Minkowski theorem.
\begin{Theorem}[Hasse-Minkowski]
 Let $A = \diag(a_1,\dots,a_n)$ and $B = \diag(b_1,\dots,b_n)$ be rational invertible matrices.
 Then $A$ and $B$ are congruent over $\mathbb Q$ (i.e.\ $B = P A P^\top$ for some $P \in \GL(n,\mathbb Q)$) if and only if the following conditions hold:
 \begin{enumerate}
  \item \label{HM:1} $A$ and $B$ have the same number of positive entries,
  \item \label{HM:2} the discriminants $\prod_{i=1}^n a_i$ and $\prod_{i=1}^n b_i$ belong to the same coset of $\mathbb Q^\times / (\mathbb Q^\times)^2$,
  \item \label{HM:3} for each odd prime $p$,
  \[
   \prod_{i<j} (a_i,a_j)_p = \prod_{i < j} (b_i,b_j)_p.
  \]
 \end{enumerate}
\end{Theorem}
If a regular antipodal plane of order $s = r-1$ exists, then $W I_v W^\top = r _v I_v$, so $I_v$ and $r I_v$ must be congruent rational matrices.
If we apply the above version of the Hasse-Minkowski theorem with $A = I_v$ and $B = r I_v$, then (\ref{HM:1}) trivially holds, (\ref{HM:2}) requires that $r^v$ is square, which holds because either $r$ is square or $v$ is even, and (\ref{HM:3}) requires that $(r,r)_p^{\binom v 2} = 1$ for all odd primes $p$.
It is easy to see from the definition of the Hilbert symbol that $(r,r)_p = (-1)^{\alpha^2 \frac{p-1}2}$ with $\alpha$ the $p$-adic valuation of $r$.
Thus, the only new condition on $r$ that we find, happens when $\binom v2$ is odd -- that is $r$ equivalent to $2$ or $4$ modulo $8$ (since we already excluded $5$ and $7$) -- and $p \equiv 3$ modulo $4$.
In that case, $p$ must divide $r$ an even number of times.
If $r \equiv 4$ modulo $8$, then we already know that $r$ must be square, so this trivially holds.
Thus, we obtain the following theorem.
\begin{Theorem}
 If $r \equiv 2$ modulo $8$, and a regular antipodal plane of order $r-1$ exists, then the square-free part of $r$ has no prime divisors congruent to $3$ modulo $4$.
\end{Theorem}

If $s\equiv1\pmod8$, then $v=(s^2+s+2)/2\equiv2\pmod4$.
Taking $r=s+1$, the preceding theorem shows that no prime congruent to
$3\pmod8$ divides the square-free part of $s+1$.
This recovers \cite[Theorem 5.2]{Payne:1972aa}.

We summarise the above results, including the Bruck-Ryser-Chowla Theorem applied to the biplane arising from a regular antipodal plane of order $s$, in the theorem below. Together with the known examples for $s=2,3$, it summarises the classical arithmetic restrictions on the orders of antipodal planes. We point out that $s=5$ was ruled out in \cite{Schellenberg1975}, and \emph{cyclic}\footnote{This means there exists a cyclic group of automorphisms acting regularly on points. We point out that \cite[p. 88]{Mielants} has misunderstood the result in \cite[p. 280]{Payne:1972aa}; only the cyclic antipodal planes for $s\le 401$ have been resolved, not all regular antipodal planes.} regular antipodal planes do not exist for $4\le s\le 401$ (see \cite[p. 280]{Payne:1972aa}). Indeed, a cyclic antipodal plane is equivalent to a \emph{$(k^2-k+2,k)$-modular Golomb ruler} where $k=s+1$ (see \cite[Section III]{Gordon}), and Gordon states that they do not exist for $k\le 10^{10}$.

\begin{Theorem}\label{badorders}
Let $s$ be the order of a regular antipodal plane. Then one of the following holds:
\begin{enumerate}[(i)]
\item $s \equiv 1 \pmod{8}$, $s-1$ is a square, no prime congruent to $3 \pmod{\textcolor{purple}{4}}$ divides the square-free part of $s+1$;
\item $s \equiv 2,5 \pmod{8}$ and $s-1$ is a square, and $s \ne 5$;
\item $s \equiv 0,7 \pmod{8}$ and $s+1$ is a square;
\item $s \equiv 3 \pmod{8}$ and $s+1$ is a square and no prime congruent to $7 \pmod{8}$ divides the square-free part of $s-1$.
\end{enumerate}
\end{Theorem}

In the next section, we will give new conditions that rule out certain orders, and give some information on the substructure of biplanes that give rise to antipodal planes.

\subsection{Biplanes}

Recall that a biplane is a symmetric $2-(v,k,2)$ design. There are only finitely many parameters for which we know that a biplane exists, and there are only 17 examples (up to isomorphism) known, including the complete $2-(4,3,2)$ design, and they have been classified for $k\le 11$ \cite{KaskiOstergard}. The known biplanes
have $k\in \{3,4,5,6,9,11,13\}$. For the purposes of regular antipodal planes, the possible known biplanes arising can only have $k\in \{3,4,6,9,11\}$ (see Corollary \ref{corollary:newcondition}). The cases $k=3,4$ exist, and we can rule out $k=5$ a number of ways. 

A \emph{signing} of a biplane $\mathcal{B}$ is a function $\sigma\colon \mathrm{Flag}(\mathcal{B})\to\{\pm 1\}$, extended by $\sigma(L,Q)=0$ for non-flags, yielding a matrix $W_\sigma\in\{-1,0,1\}^{\mathcal{L}\times\mathcal{P}}$. The signing is a \emph{weighing signing} if $W_\sigma (W_\sigma)^\top = kI_v$.

Moorhouse \cite{Moorhouse} studied the reconstruction of projective planes from homology semibiplanes, reducing the choices in the double cover of the nonfixed part to linear equations over $\mathbb{F}_2$. The lemma below applies the same formulation to biplanes.

\begin{Lemma}\label{lem:f2}
Write $\sigma = (-1)^x$ for $x\colon \mathrm{Flag}(\mathcal{B})\to\mathbb{F}_2$. Then $\sigma$ is a weighing signing if and only if $x$ satisfies, for every unordered pair $\{L_1,L_2\}$ of distinct lines, meeting in two points $P_1$ and $P_2$,
\begin{equation}\label{eq:sys}
x(L_1,P_1) + x(L_2,P_1) + x(L_1,P_2) + x(L_2,P_2) \equiv 1 \pmod 2.
\end{equation}
\end{Lemma}

\begin{proof}
Each row of $W_\sigma$ has exactly $k$ nonzero entries, all $\pm 1$, so each diagonal entry of $(W_\sigma(W_\sigma)^\top)$ is equal to $k$.
For distinct lines $L_1,L_2$, meeting in points $P_1$ and $P_2$,
\[
(W_\sigma(W_\sigma)^\top)[L_1,L_2] = \sum_{Q\in L_1\cap L_2}\sigma(L_1,Q)\sigma(L_2,Q) = \sigma(L_1,P_1)\sigma(L_2,P_1) + \sigma(L_1,P_2)\sigma(L_2,P_2).
\]
Both summands lie in $\{\pm 1\}$, so the sum vanishes if and only if the summands are opposite signs, that is,
\[
\sigma(L_1,P_1)\sigma(L_2,P_1)\sigma(L_1,P_2)\sigma(L_2,P_2) = -1.
\]
Substituting $\sigma=(-1)^x$ gives exactly~\eqref{eq:sys}.
\end{proof}

It is known via a result of \cite{Tolonen} that there is no regular antipodal plane of order 8, but only once we establish a connection to equidistant permutation arrays in the next section. We show that we can use biplanes to establish this result.

\begin{Theorem}\label{8or10}
    There is no regular antipodal plane of order 8 or 10.
\end{Theorem}

\begin{proof}
    If such a regular antipodal plane exists, then its antipodal quotient is a biplane
    of order $7$ or $9$ admitting a weighing signing.  Let $\mathcal B$
    be such a biplane, and let $\Gamma$ be its bipartite incidence graph.
    By Lemma~\ref{lem:f2}, a weighing signing is equivalent to a function
    $x\colon \mathrm{Flag}(\mathcal B)\to \mathbb F_2$
    such that, for every pair of distinct lines $L_1,L_2$, meeting in the
    two points $P_1,P_2$,
    \[
        x(L_1,P_1)+x(L_2,P_1)+x(L_1,P_2)+x(L_2,P_2) \equiv 1 \pmod 2.
    \]
    Thus we need to decide whether this affine system over $\mathbb F_2$
    is soluble.

    Take the flags of $\mathcal B$, equivalently the edges $E$ of
    $\Gamma$, together with a new symbol $\infty$.  Let $W$ be the
    $\mathbb F_2$-vector space with basis $E\cup\{\infty\}$.  For each
    pair of distinct lines $L_1,L_2$, meeting in $P_1,P_2$, let
    $e_{ij}$ denote the edge $(L_i,P_j)$ of $\Gamma$, and form the
    vector
    \[
        e_{11}+e_{21}+e_{12}+e_{22}+\infty .
    \]
    Let $U$ be the subspace of $W$ generated by all such vectors.

    The affine system above is soluble if and only if $\infty\notin U$.
    Indeed, a solution $x$ defines a linear functional $\varphi\colon W\to
    \mathbb F_2$ by
    \[
        \varphi(e)=x(e)\quad(e\in E),\qquad \varphi(\infty)=1,
    \]
    and the equations say exactly that $\varphi$ vanishes on the generators
    of $U$.  Hence $\varphi$ vanishes on $U$, so $\infty\notin U$.
    Conversely, if $\infty\notin U$, then there is a linear functional
    $\varphi\colon W\to\mathbb F_2$ vanishing on $U$ but with
    $\varphi(\infty)=1$; restricting $\varphi$ to $E$ gives a solution
    of the system.

    For each biplane of order $7$ and each biplane of order $9$, we
    constructed the corresponding subspace $U$ and row-reduced over
    $\mathbb F_2$.  In every case the vector $\infty$ lies in $U$.
    Therefore the system from Lemma~\ref{lem:f2} has no solution, and hence
    no such biplane admits a weighing signing.  Consequently no regular
    antipodal plane of order $8$ or $10$ exists.
\end{proof}

If a regular antipodal plane of order $s$ exists, where
$2\le s\le 14$, then Theorem \ref{badorders} leaves only $2,3,8,10$ in that interval; Theorem \ref{8or10} removes $8,10$. So we have:

\begin{Corollary}
    A regular antipodal plane of order $2\le s\le 14$ exists if and only
    if $s\in\{2,3\}$.
\end{Corollary}

\begin{Question}
    Does every finite $K_{2,3}$-free bipartite graph occur as a subgraph of the incidence graph of some finite biplane?
\end{Question}

All subgraphs of the incidence graph of a biplane are $K_{2,3}$-free. As of the writing of this paper, there is no finite $K_{2,3}$-free bipartite graph that has been shown to not be a subgraph of any incidence graph of a finite biplane.

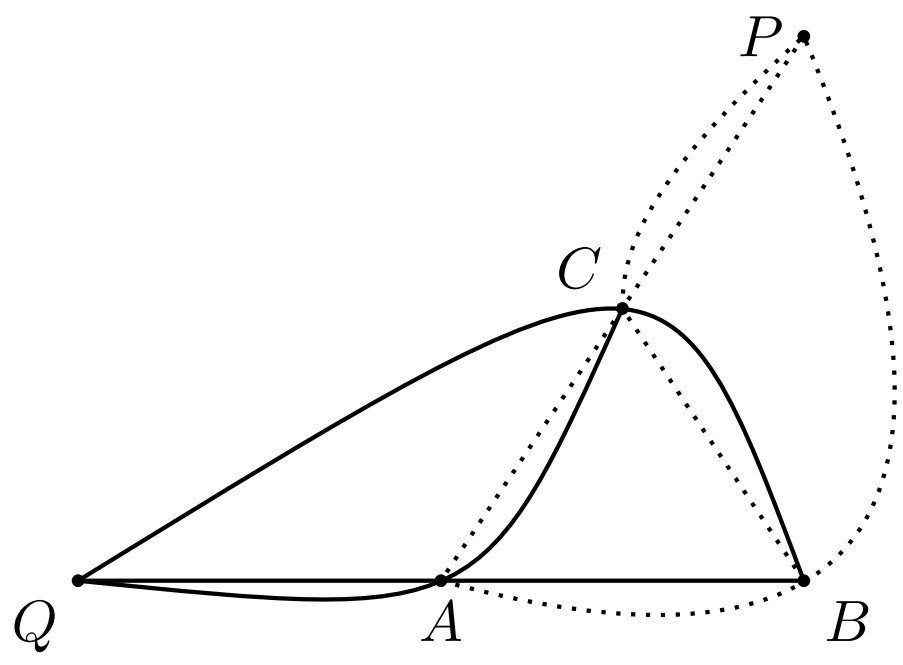
\begin{figure}[!ht]
\centering
\begin{tikzpicture}[scale=0.6]
  \coordinate (A) at (0,0);
  \coordinate (B) at (4,0);
  \coordinate (C) at (2,3);
  \coordinate (P) at (4,6);
  \coordinate (Q) at (-4,0);

  \draw[dotted, thick]
    (A) .. controls (4,-1.1) and (6.5,0) .. (P);

  \draw[dotted, thick]
    (C) .. controls (2.0,3.5) and (2.1,4.3) .. (P);

  \draw[thick]
    (B) .. controls (2.5,4) .. (Q);

  \draw[thick]
    (C) .. controls (.45,-.45) .. (Q);


   \draw[ thick] (-4,0) -- (4,0);
  \draw[dotted, thick] (4,0) -- (2,3);
  \draw[dotted, thick] (4,6) -- (0,0);

  \fill (P) circle (2pt);
    \fill (A) circle (2pt);
    \fill (Q) circle (2pt);
  \fill (B) circle (2pt);
  \fill (C) circle (2pt);
\node[below] at (A) {$A$};
\node[below right] at (B) {$B$};
\node[above left] at (C) {$C$};
\node[left] at (P) {$P$};
\node[below left] at (Q) {$Q$};

\end{tikzpicture}
\caption{A minimal forbidden $K_{2,3}$-free graph on 11 vertices: 5 points and 6 lines.}\label{fig:5points6lines}
\end{figure}

However, for a biplane that is the quotient of an antipodal plane, the linear system constraint does forbid many subgraphs. Below we give all of the minimal forbidden $K_{2,3}$-free graphs on at most 16 vertices, as well as a picture of the unique smallest one on 11 vertices, with 5 points and 6 lines (Figure \ref{fig:5points6lines}). For brevity, they are given as sets of lines of an incidence structure with $v$ points. Note that from the point of view of incidence geometry, both the given configurations and their duals are forbidden. 

The known biplanes up to duality in Table \ref{biplanes} are denoted by $bi$ where $i = k-2$, also known as the order of the biplane.  When there is more than one nonisomorphic example, the letters $A,B,C,\dots$ are appended to distinguish them. The correspond in the literature as follows: $b4_A,b4_B,b4_C$ correspond to $B_6, B_7, B_8$ in \cite{AssmusSalwach}, $b7_A,b7_B,b7_C$
correspond to $B_H(9), B_L(9), B_d(9)$ in \cite{Salwach}, respectively.  The biplanes $b9_A,b9_B,b9_C,b9_D,b9_E$ correspond to $\mathcal{B}_1, \mathcal{B}_2,\mathcal{B}_4,\mathcal{B}_3,\mathcal{B}_5$ in \cite{KaskiOstergard}, respectively.
All the planes except order 11 can be downloaded from \cite{Heinlein}.  
The known biplane of order $11$ can be found in \cite{Aschbacher}. 

\begin{table}[!ht]\label{biplanes}
{\footnotesize
\begin{tabular}{c||c|c|c|c|c|c|c|c|c|c|c|c|c|c|c|c|c|c}
\toprule
& $v$ & Aut. group & $C_1$ & $C_2$ & $C_3$ & $C_4$ & $C_5$ & $C_6$ & $C_7$ & $C_8$ & $C_9$ & $C_{10}$ & $C_{11}$ & $C_{12}$ & $C_{13}$ & $C_{14}$ & $C_{15}$ & $C_{16}$ \\
\midrule
$b1$ & 4 &  $S_4$ &  &   &   &   &   &   &   &   &   &   &   &   &   &   &   &   \\
$b2$ & 7 &  $\mathrm{PSL}(3,2)$ & &   &   &   &   &   &   &   &   &   &   &   &   &   &   &   \\
$b3$ & 11 & $\mathrm{PSL}(2,11)$ & $\times$ &   &   &  $\times$ &   &   &   &   &   &   &   &   &   &   &   &  $\times$ \\
$b4_A$ & 16 & $2^4:S_6$ & &  &   &   &   &   &   &   &   &   &   &   &   &   &   &  $\times$ \\
$b4_B$ & 16 & $2^5\cdotp S_4$&  & &  $\times$ &  $\times$ &   &   &   &   &   &   &   &   &  $\times$ &   &   &  $\times$ \\
$b4_C$ & 16 & $2^{1+4}_-:A_4$& &  &  $\times$ &  $\times$ &   &   &  $\times$ &   &   &   &   &   &  $\times$ &   &   &   \\
$b7_A$ & 37 & $\mathrm{P\Gamma L}(2,8)$ & $\times$ &  $\times$ &  $\times$ &   &   &   &  $\times$ &  $\times$ &   &  $\times$ &  $\times$ &  $\times$ &   &  $\times$ &  $\times$ &   \\
$b7_B$ & 37 & $C_3^2: C_6$ & $\times$ &  $\times$ &  $\times$ &   &  $\times$ &  $\times$ &  $\times$ &  $\times$ &  $\times$ &  $\times$ &  $\times$ &  $\times$ &  $\times$ &  $\times$ &  $\times$ &   \\
$b7_C$ & 37 & $C_{37} : C_9$ & &  $\times$ &   &   &  $\times$ &  $\times$ &  $\times$ &  $\times$ &  $\times$ &  $\times$ &   &  $\times$ &  $\times$ &  $\times$ &   &   \\
$b9_A$ & 56 & $\mathrm{PSL}(3,4) : C_2^2$ & &   &   &   &   &   &   &   &   &   &   &   &   &   &   &  $\times$ \\
$b9_B$ & 56 & $C_2 \times (C_3^2 : SD_{16})$ & $\times$ &  $\times$ &  $\times$ &  $\times$ &   &  $\times$ &  $\times$ &  $\times$ &  $\times$ &  $\times$ &  $\times$ &  $\times$ &  $\times$ &  $\times$ &  $\times$ &  $\times$ \\
$b9_C$ & 56 & $(C_2 \times SD_{16}) : C_2$ & $\times$ &  $\times$ &  $\times$ &  $\times$ &  $\times$ &  $\times$ &  $\times$ &  $\times$ &  $\times$ &  $\times$ &  $\times$ &  $\times$ &  $\times$ &  $\times$ &  $\times$ &  $\times$ \\
$b9_D$ & 56 & $C_3^2: SD_{16}$ & &  $\times$ &  $\times$ &   &  $\times$ &  $\times$ &  $\times$ &  $\times$ &  $\times$ &  $\times$ &  $\times$ &  $\times$ &  $\times$ &  $\times$ &  $\times$ &  $\times$ \\
$b9_E$ & 56 & $C_2 \times A_4$ & $\times$ &  $\times$ &  $\times$ &  $\times$ &  $\times$ &  $\times$ &  $\times$ &  $\times$ &  $\times$ &  $\times$ &  $\times$ &  $\times$ &  $\times$ &  $\times$ &  $\times$ &  $\times$ \\
$b11$ & 79 & $\mathrm{AGL}(1,11)$ & $\times$ &  $\times$ &  $\times$ &   &  $\times$ &  $\times$ &  $\times$ &  $\times$ &  $\times$ &  $\times$ &  $\times$ &  $\times$ &  $\times$ &  $\times$ &  $\times$ &   \\
\bottomrule
\end{tabular}}\smallskip

\caption{Each row denotes a known biplane, up to duality. 
Each $C_i$ is a possible minimal obstruction to a signing, and a $\times$ denotes an occurrence in the
relevant biplane. We have listed the automorphism groups of the biplanes here, for easier identification.}
\end{table}

\begin{table}[!ht]
\begin{tabular}{l||l}
\toprule
     $C_{1}$&  $\{ 1, 2, 4 \}, \{ 1, 2, 5 \}, \{ 1, 3, 4 \}, \{ 1, 3, 5 \}, \{ 2, 4, 5 \}, \{ 3, 4, 5 \}$\\
     $C_{2}$&  $\{ 1, 2, 5 \}, \{ 1, 3, 4 \}, \{ 1, 3, 5 \}, \{ 2, 3, 5 \}, \{ 2, 3, 6 \}, \{ 3, 4, 6 \}, \{ 1, 2, 4, 6\}$\\
     $C_{3}$&  $\{ 1, 2, 5 \}, \{ 1, 2, 6 \}, \{ 1, 3, 5 \}, \{ 1, 3, 6 \}, \{ 2, 4, 5 \}, \{ 2, 4, 6 \}, \{ 3, 4, 5, 6 \}$\\
     $C_{4}$&  $\{ 1, 2, 4 \}, \{ 1, 3, 4 \}, \{ 4, 5, 6 \}, \{ 4, 5, 7 \}, \{ 1, 2, 5, 6 \}, \{ 1, 3, 5, 7 \}, \{ 2, 3, 4, 6, 7 \}$\\
     $C_{5}$&  $\{ 1, 4, 5 \}, \{ 1, 4, 6 \}, \{ 2, 3, 4 \}, \{ 2, 3, 5 \}, \{ 2, 4, 5 \}, \{ 3, 5, 7 \}, \{ 1, 5, 6, 7 \}, \{ 3, 4, 6, 7 \}$\\
     $C_{6}$&  $\{ 1, 3, 4 \}, \{ 1, 4, 5 \}, \{ 2, 3, 4 \}, \{ 2, 4, 5 \}, \{ 2, 5, 7 \}, \{ 5, 6, 7 \}, \{ 1, 3, 5, 6 \}, \{ 2, 3, 6, 7 \}$\\
     $C_{7}$&  $\{ 1, 3, 4 \}, \{ 1, 3, 5 \}, \{ 1, 4, 5 \}, \{ 2, 3, 5 \}, \{ 2, 5, 7 \}, \{ 2, 6, 7 \}, \{ 2, 3, 4, 6 \}, \{ 4, 5, 6, 7 \}$\\
     $C_{8}$&  $\{ 1, 2, 5 \}, \{ 1, 3, 5 \}, \{ 2, 4, 5 \}, \{ 3, 4, 6 \}, \{ 4, 5, 7 \}, \{ 4, 6, 7 \}, \{ 1, 2, 3, 4 \}, \{ 3, 5, 6, 7 \}$\\
     $C_{9}$&  $\{ 1, 2, 5 \}, \{ 1, 4, 5 \}, \{ 2, 5, 7 \}, \{ 3, 4, 5 \}, \{ 3, 4, 6 \}, \{ 3, 5, 7 \}, \{ 1, 2, 4, 6 \}, \{ 2, 3, 6, 7 \}$\\
     $C_{10}$&  $\{ 1, 2, 5 \}, \{ 1, 3, 5 \}, \{ 2, 4, 6 \}, \{ 2, 5, 6 \}, \{ 3, 4, 7 \}, \{ 3, 5, 7 \}, \{ 1, 2, 3, 4 \}, \{ 4, 5, 6, 7 \}$\\
     $C_{11}$&  $\{ 1, 2, 3 \}, \{ 1, 2, 4 \}, \{ 1, 3, 4 \}, \{ 2, 4, 6 \}, \{ 3, 4, 7 \}, \{ 4, 5, 6 \}, \{ 4, 5, 7 \}, \{ 2, 3, 5, 6, 7 \}$\\
     $C_{12}$&  $\{ 1, 3, 4 \}, \{ 1, 3, 7 \}, \{ 1, 4, 5 \}, \{ 2, 3, 4 \}, \{ 2, 3, 7 \}, \{ 2, 4, 6 \}, \{ 4, 5, 6 \}, \{ 1, 2, 5, 6, 7 \}$\\
     $C_{13}$&  $\{ 1, 3, 4 \}, \{ 1, 3, 6 \}, \{ 1, 4, 7 \}, \{ 2, 3, 4 \}, \{ 2, 4, 7 \}, \{ 2, 5, 7 \}, \{ 1, 5, 6, 7 \}, \{ 2, 3, 5, 6 \}$\\
     $C_{14}$&  $\{ 1, 3, 5 \}, \{ 1, 6, 7 \}, \{ 2, 3, 4 \}, \{ 2, 3, 5 \}, \{ 2, 4, 7 \}, \{ 4, 6, 7 \}, \{ 1, 2, 5, 7 \}, \{ 1, 3, 4, 6 \}$\\
     $C_{15}$&  $\{ 1, 2, 5 \}, \{ 1, 3, 5 \}, \{ 2, 4, 5 \}, \{ 3, 5, 6 \}, \{ 4, 5, 7 \}, \{ 5, 6, 7 \}, \{ 1, 2, 3, 4 \}, \{ 3, 4, 6, 7 \}$\\
    $C_{16}$&  $\{ 1, 2, 3 \}, \{ 1, 2, 6 \}, \{ 1, 3, 5 \}, \{ 1, 4, 5 \}, \{ 1, 4, 6 \}, \{ 2, 3, 4 \}, \{ 2, 4, 5 \}, \{ 2, 5, 6 \}, \{ 3, 4, 6 \}, \{ 3, 5, 6 \}$\\
 \bottomrule
 \end{tabular}\smallskip
 \caption{Forbidden configurations.}
\end{table}

This list was found in two ways, by first using Brendan McKay's graph generation program \texttt{geng} \cite{nauty} and independently by examining the subconfigurations of the known biplanes using GAP. Both lists were found to be equal.  In each configuration,  the system described in Lemma \ref{lem:f2} is inconsistent, in particular the sum of all of the corresponding equations equals 1.

\begin{Proposition}\label{propn:oddorder}
    Let $\mathcal{B}$ be a biplane. If $\mathrm{Aut}(\mathcal{B})$ contains a subgroup $G$ of odd order acting transitively on $\mathrm{Flag}(\mathcal{B})$, then $\mathcal{B}$ admits no weighing signing.
\end{Proposition}

\begin{proof}
 Consider the $\mathbb F_2$-matrix $M$ of the linear system given by (\ref{eq:sys}), that is the columns of $M$ are indexed by the flags of $\mathcal B$, and for any pair of lines $\{L_1,L_2\}$ with common points $\{P_1,P_2\}$, $M$ has a row with ones exactly in positions $(P_1,L_1), (P_1,L_2), (P_2,L_1), (P_2,L_2)$.
 A weighing signing of $\mathcal B$ corresponds to a vector $x$ with $M x = \boldsymbol 1$.
 Note that $G$ naturally acts on $\mathbb F_2^{{\rm Flag}(\mathcal B)}$, and $M x^g = \boldsymbol 1$ for each $g \in G$.
 Define $\overline x = \sum_{g \in G} x^g$.
 Then $M \overline x = |G| \boldsymbol 1 = \boldsymbol 1$.
 On the other hand, by the transitivity of $G$ on ${\rm Flag}(\mathcal B)$, $\overline x$ must be constant, that is $\overline x = a \boldsymbol 1$ for some $a \in \mathbb F_2$.
 But then $M \overline x = a M \boldsymbol 1 = \boldsymbol 0$, because every row of $M$ has four ones.
 This yields a contradiction.
\end{proof}

The only known biplanes satisfying the hypothesis of Proposition \ref{propn:oddorder} are (i) the unique biplane on 11 points ($k=5$, $s=4$)
(admitting $C_{11}:C_5\le \PSL(2,11)$), and (ii) a biplane on 37 points ($k=9$, $s=8$) admitting $C_{37}:C_9$. 
Another way to constrain the construction of future examples is to use information about the $\F_2$-rank of
a putative example:

\begin{Proposition}\label{exact-rank}
Let $\mathcal D$ be a biplane admitting a weighing signing, with block size $k$, and let
$v=1+\binom{k}{2}$. If $k$ is even, then the incidence matrix $B$ of $\mathcal D$ satisfies
\begin{equation}
  \rank_{\F_2}B=\left\lfloor\frac v2\right\rfloor.
  \label{eq:exact-rank}
\end{equation}
\end{Proposition}

\begin{proof}
Let $W$ be a weighing signing of $B$, and let $N$ be the negative support of
$W$; thus
\[
  N_{ij}=1\iff W_{ij}=-1,
  \qquad W=B-2N.
\]
Expanding $WW^\top=kI$ and using $BB^\top=B^\top B=(k-2)I+2J$, we obtain 
\[
  BN^\top+NB^\top-2NN^\top=J-I.
\]
Reducing modulo $2$ gives
$BN^\top+NB^\top=I+J$. 
All ranks are over $\F_2$. Let $r:=\rank B$. Then
\[
2\left\lfloor\frac v2\right\rfloor=\rank(I+J)\le \rank(BN^\top)+\rank(NB^\top)\le 2r.
\]
Since $k$ is even, the relation $BB^\top=B^\top B=(k-2)I+2J$ gives
$B^\top B=BB^\top=0$. So the row-space of $B$ is self-orthogonal and hence $r\le \left\lfloor\frac v2\right\rfloor$,
and equality follows.
\end{proof}

Proposition \ref{exact-rank} gives us another method to rule out certain regular antipodal planes,
when $k$ is even. For $k=6$ there are exactly three biplanes. Assmus and Salwach \cite{AssmusSalwach} denote them
according to the dimensions of their binary incidence codes; we write
$\mathcal D_6,\mathcal D_7,\mathcal D_8$, with
$\rank_2 B(\mathcal D_i)=i$. Proposition~\ref{exact-rank} requires rank $v/2=8$, and so
$\mathcal{D}_6$ and $\mathcal{D}_7$ are not signable.
Using standard representatives of the classified biplanes gives the
following exact ranks.

\begin{table}[ht]
\centering
\label{tab:small-tests}
\begin{tabular}{@{}lrrcc@{}}
\toprule
Design & $k$ & $v$ & $\rank_2 B$ & Required rank for a weighing signing \\
\midrule
$\mathcal D_4$ & 4 & 7  & 3  & 3  \\
$\mathcal D_6$ & 6 & 16 & 6  & 8  \\
$\mathcal D_7$ & 6 & 16 & 7  & 8  \\
$\mathcal D_8$ & 6 & 16 & 8  & 8  \\
\bottomrule
\end{tabular}\smallskip
\caption{Using Proposition~\ref{exact-rank} on some known biplanes.}
\end{table}

Thus the exact-rank theorem eliminates two of the three $2-(16,6,2)$
biplanes without solving any signing equations. No other known biplane is
eliminated with this method.

\subsection{Sharply deranged sets of permutations}\label{sharplyderanged}

For integers $a \le b$, denote $\{a,\dots,b\}$ by $[a,b]$.
Schleiermacher \cite[Corollary 5.3]{Schleiermacher:2016aa} showed that
there exists a finite regular antipodal plane of order $s\ge 2$ if and only if there exists a set $\Phi$ of permutations in $S_{s+1}$ satisfying: (i) every element of $\Phi$ has no cycles of length at most 2; (ii) if $(a, b)$ and $(c,d)$ are two ordered pairs of elements of $[1,s+1]$ such that $a\ne b$, $c\ne d$, $a\ne c$, $b\ne d$, and $(c,d)\ne (b,a)$, then there exists a unique element $g\in \Phi$ such that $b=a^g$ and $d=c^g$.
We show here that there is a simpler set of permutations that is equivalent to a regular antipodal plane.
We call a family $\mf$ of permutations in $S_{s+1}$ \emph{sharply deranged} if
\begin{enumerate}[(i)]
 \item $|\mf| = s(s-1)$,
 \item no element of $\mf$ stabilises a point or pair setwise,
 \item for every $g \in \mf$, it holds that $g^{-1} \in \mf$, and for every $f \in \mf \setminus \{g,g^{-1}\}$, $gf^{-1}$ has a unique fixed point.
\end{enumerate}

\begin{Proposition}
 A regular antipodal plane of order $s \ge 2$ exists if and only if $S_{s+1}$ has a sharply deranged family.
\end{Proposition}

\begin{proof}
Suppose that a regular antipodal plane of order $s$ exists. Take two points $P$ and $P^*$. Let $\ell_0, \dots, \ell_s$ be the lines incident with $P$. Then every point $Q \ne P, P^*$ is on a unique line $\ell_i$ and a unique line $\ell_j^*$. Moreover, $i \ne j$ since $\ell_i$ and $\ell_i^*$ don't intersect. Vice versa, if $i \ne j$, then $\ell_i$ and $\ell_j^*$ intersect in a unique point. We denote this point by $(i,j)$. Thus, the set of points of the plane is 
 \[
  ( [0,s]^2 \setminus \{(i,i) : i \in [0, s] \} ) \cup \{P,P^*\}.
 \] 
Note that $(\ell_i \cap \ell_j^*)^* = \ell_j \cap \ell_i^*$, hence $(i,j)^* = (j,i)$.

Now take a line $\ell$ of the plane, not incident with $P$ or $P^*$.
Then $\ell$ intersects each line $\ell_i$ and $\ell_i^*$ in exactly one point, hence the points of $\ell$ are of the form $\{(i,g_\ell(i)) : i \in [0,s]\}$ for a permutation $g_\ell$. Note that $g_\ell$ has no fixed points since $(i,i)$ is not a point of the plane, and $g_\ell$ does not stabilise a pair $\{i,j\}$, since $(i,j)$ and $(j,i) = (i,j)^*$ cannot both be incident with $\ell$. Furthermore, $(i,j)$ is incident with $\ell$ if and only if $(i,j)^* = (j,i)$ is incident with $\ell^*$, hence, $g_{\ell^*} = g_\ell^{-1}$. The family of permutations $g_\ell$ on $[0,s]$ is sharply deranged.

Conversely, suppose that $S_{s+1}$ is the symmetric group acting on $[0,s]$, and that $\mf$ is a sharply deranged family.
Make an incidence structure with point set
\[
\mp = \{P,P^*\} \cup \{ (i,j) \in [0,s]^2 : i \ne j \},
\]
and line set 
\[
\ml = \{\ell_i : i \in [0,s] \} \cup \{ \ell_i^* : i \in [0,s]\} \cup \mf.
\]
The points incident with line $\ell$ are
\begin{enumerate}[(i)]
\item if $\ell = \ell_i$: $P$ and the points $(i,j)$ with $j \in [0,s] \setminus \{i\}$, 
\item if $\ell = \ell_i^*$: $P^*$ and the points $(j,i)$ with $j \in [0,s] \setminus \{i\}$, 
\item if $\ell = g \in \mf$: the points $(i,g(i))$ with $i \in [0,s]$.
\end{enumerate}
It follows that every line is incident with $s+1$ points, and two distinct lines intersect in at most one point.
We show that every point is incident with $s+1$ lines.
There are equally many lines as points, and every line is incident with $s+1$ points, so on average a point is incident with $s+1$ lines.
Now suppose that some point $Q$ is incident with more than $s+1$ lines.
Then these lines pairwise intersect in $Q$, hence there are at least $1 + (s+2)s > s^2 + s + 2$ points, which is a contradiction.
Therefore, every point must be incident with $s+1$ lines.
This proves that we have constructed a regular antipodal plane.
\end{proof}

\begin{Example}
    For $s=2$, $\{(1,2,3),(1,3,2)\}$ is sharply deranged. Likewise, the set of 4-cycles is sharply deranged for $s=3$.
    For $s=4$, you cannot extend a 5-cycle plus its inverse to anything larger. 
\end{Example}

Suppose $\mf$ is a sharply deranged family of permutations in $S_{s+1}$. For each inverse pair $(g,g^{-1})$ in $\mf$, throw out one of them, which will in total give us a set $\mf'$ of half the elements of $\mf$. Now each pair of distinct permutations agree on precisely one point, and hence $\mf'$ is an \emph{equidistant permutation array} (EPA).  Table \ref{tbl:comparison} compares the clique number of a graph $\Omega_{s+1}$ which we will define below, and the size $R(s+1,1)$ of the largest EPA in $S_{s+1}$ (see also OEIS:A005677). In particular, Tolonen \cite{Tolonen} showed recently that $R(9,1) \le 20$, from which it follows that there is no sharply deranged family of permutations in $S_9$, and hence no regular antipodal plane of order 8. The graph $\Omega_{s+1}$ has vertices the derangements of $S_{s+1}$ that do not have 2-cycles. Two permutations $g,h$ are adjacent in $\Omega_{s+1}$ if $g=h^{-1}$ or $gh^{-1}$ fixes a unique element. We computed clique numbers by computer, using the GAP \cite{GAP} package \textsc{Grape} \cite{GRAPE}.

\begin{table}[!ht]
\begin{tabular}{lllll}
\toprule
  $s$  & $s(s-1)$ & $\omega(\Omega_{s+1})$ & $s(s-1)/2$ & $R(s+1,1)$  \\
\midrule
   3  & 6 & 6 & 3 & 3\\    
   4 & 12 & 2 & 6 & 5\\
   5 & 20 & 11 & 10 & 10\\
   6 & 30 & 10 & 15 & 13\\
   7 & 42 & 14 & 21 &16\\
   8 & 56 & $\ge 20$ & 28 &20\\
\bottomrule\\[-4pt]
\end{tabular}
\caption{Clique numbers and bounds. If a sharply deranged subset of $S_{s+1}$ exists, then $s(s-1)\le \omega(\Omega_{s+1})$
and $s(s-1)/2\le R(s+1,1)$.}\label{tbl:comparison}
\end{table}

The data in \cite{Mahmoodi:1996aa} shows that the known EPAs would easily satisfy an upper bound of $\tfrac{1}{4}n^2+2$. So based on the numerical evidence, we propose the following conjecture:

\begin{Conjecture}\label{EPAconjecture}
    $R(n,1) \le \frac{1}{4}n^2+2$ for all $n\ge 3$.
\end{Conjecture}

The best bound on an EPA is by van Rees \cite{Van-Rees:1990aa}: $n^2-4n+O(\sqrt n)$. It was shown in \cite{MathonVanstone} that $R(n,1)\ge \frac{1}{2}(n-1)(\sqrt{4n-3}-1)$ when $n=q^2+q+1$ for some prime-power $q$.  If Conjecture \ref{EPAconjecture} holds, then so too does Payne's Conjecture.
Choosing one element from each inverse pair in a sharply deranged set of $S_{s+1}$ yields an EPA of index 1 of size $s(s-1)/2$. If Conjecture \ref{EPAconjecture} holds, then $\frac{1}{2} s(s-1)\le \frac{1}{4}(s+1)^2+2$, from which it follows that $s\le 5$. (We have already excluded the cases $s=4,5$).

\subsection{Mutually orthogonal array-bic squares}\label{arraybic}

\begin{Definition}
 An \emph{array-bic square} of order $s$ and position $i$, with $1 \le i \le s$ is an $s \times s$-array $A_i$ with the following properties:
 \begin{enumerate}[(i)]
  \item The diagonal of $A_i$ is empty, the off-diagonal positions in $A_i$ contain entries from $[0,s]\setminus\{i\}$.
  \item The $i$th column of $A_i$ has a $0$ entry on all off-diagonal positions,
  \item The $y$th column of $A_i$, with $y \neq i$, contains each symbol in $[1,s] \setminus \{i\}$ exactly once.
  \item The $x$th row of $A_i$ contains each symbols of $[1,s]\setminus\{i,x\}$ exactly once, and if $i \neq x$, it also contains a $0$.
 \end{enumerate}
\end{Definition}

Here is an example of an array-bic square of order 4 and position 1:
\[
\begin{bmatrix}
\cdot&2&4&3\\
0&\cdot&3&4\\
0&4&\cdot&2\\
0&3&2&\cdot\\
\end{bmatrix}
\]

\begin{Definition}
 Suppose that $A_i$ and $A_k$ are array-bic squares of order $s$ and positions $i$ and $k$, respectively, with $i\ne k$.
 We call $A_i$ and $A_k$ \emph{orthogonal} if for every tuple 
$(a,b)\in ([0,s]\setminus\{i\})\times([0,s]\setminus\{k\})$ other than
\[
    (k,i),\qquad
    (c,c)\quad\text{for }c\in [0,s]\setminus\{i,k\},
\]
there is a unique off-diagonal position $(x,y)$ such
that
\[
    A_i(x,y)=a  \qquad\text{and}\qquad A_k(x,y)=b.
\]
(Equivalently, the pairs which do not occur are precisely
$(k,i)$ and $(c,c)$ for $c\in [0,s]\setminus\{i,k\}$.)
 Note that this accounts for all $s(s-1)$ off-diagonal positions.
 Note that if $i=k$, then there are positions where we find $(a,b) = (0,0)$, so $A_i$ and $A_k$ cannot be orthogonal.
 
 We call $A_{i_1}, \dots, A_{i_t}$ a set of \emph{mutually orthogonal array-bic squares}, abbreviated $\MOAS(s,t)$, if they are array-bic squares of order $s$, and any pair of them is orthogonal.
\end{Definition}

\begin{Example}
For order 2, there is a unique example of a $\MOAS(2,2)$:
\[
\begin{bmatrix}
\cdot&2\\
0&\cdot
\end{bmatrix},
\begin{bmatrix}
\cdot&0\\
1&\cdot
\end{bmatrix}
\]
By computer, it is not difficult to show that the unique $\MOAS(3,3)$ (up to symmetry) is
\[
\begin{bmatrix}
\cdot&3&2\\
0&\cdot&3\\
0&2&\cdot\\
\end{bmatrix},
\begin{bmatrix}
\cdot&0&3\\
3&\cdot&1\\
1&0&\cdot\\
\end{bmatrix},
\begin{bmatrix}
\cdot&2&0\\
1&\cdot&0\\
2&1&\cdot\\
\end{bmatrix}
\]
\end{Example}

The proof of the following is straightforward from the definitions, but it helps to know this result in
efficient computation and enumeration of examples.

\begin{Lemma}\label{Lemma:arraybic_auts}
 Suppose that $A_i$ is an array-bic square of order $s$ and position $i$.
 Consider the natural action of the symmetric group $S_s$ on $[1,s]$, 
and stipulate that it fixes the symbol $0$. Take $g \in S_s$. 
 Define $A_i^g$ by
 \[ 
  A_i^g( x, y ) = \left( A_i( x^{g^{-1}}, y^{g^{-1}} ) \right)^g.
 \]
 Then $A_i^g$ is an array-bic square of order $s$ and position $i^g$.
 Moreover, two array-bic squares $A_i$ and $A_k$ are orthogonal if and only if $A_i^g$ and $A_k^g$ are orthogonal.
\end{Lemma}

\begin{Proposition}
 A regular antipodal plane of order $s$ is equivalent to a $\MOAS(s,s)$.
\end{Proposition}

\begin{proof}
 Suppose that $(\mp, \ml)$ is a regular antipodal plane of order $s$.
 Choose a line $\ell_0 \in \ml$, and label the points of $\ell_0$ as $P_0, \dots, P_s$.
 For each $(i,j) \in [0,s]^2$ with $i \neq j$, let $\ell_{ij}$ be the line containing $P_i$ and $P_j^*$.
 Note that $\ml = \{ \ell_0, \ell_0^* \} \cup \{ \ell_{ij} : i,j \in [0,s], i \neq j \}$, and $\ell_{ij}^* = \ell_{ji}$.
 For each $(x,y) \in [1,s]^2$ with $x \neq y$, we denote the point $ \ell_{0x} \cap \ell_{y0}$ by $P_{xy}$.
 Note that $P_{xy}^* = P_{yx}$.
 Define the $s \times s$ arrays $A_1, \dots, A_s$ with empty diagonal by $A_i(x,y) = j$ if the unique line containing $P_i$ and $P_{xy}$ is $\ell_{ij}$.
 Using the properties of regular antipodal planes, it is straightforward to verify that each $A_i$ is an array-bic square.
 Suppose by way of contradiction that $A_i$ and $A_k$ are not orthogonal.
 If $A_i(x,y)=k$ and $A_k(x,y)=i$, then $P_{xy}$ lies on both
$\ell_{ik}$ and $\ell_{ki}=\ell_{ik}^*$, contradicting their
disjointness. So one of the following must occur:
 \begin{enumerate}
  \item There exists a off-diagonal position $(x,y)$ with $A_i(x,y) = A_k(x,y) = j$.
  This means that the line through $P_j^*$ and $P_{xy}$ intersects $\ell_0$ in both $P_i$ and $P_k$, a contradiction.
  \item There exists two off-diagonal positions $(x,y), (x',y')$ with $A_i(x,y) = A_i(x',y') = j$ and $A_k(x,y) = A_k(x',y') = h$ for some $j$ and $h$.
  This would mean that the lines $\ell_{ij}$ and $\ell_{kh}$ share the points $P_{xy}$ and $P_{x'y'}$, a contradiction.
 \end{enumerate}
 We conclude that a regular antipodal plane gives rise to a $\MOAS(s,s)$.

Conversely, suppose that $A_1,\ldots,A_s$ is a
$\MOAS(s,s)$, where $A_i$ has position $i$.
Put
\[
    X=[0,s],
    \qquad
    \Omega=\{(x,y)\in[1,s]^2:x\ne y\},
\]
and define maps $f_0(x,y)=x$, $f_i(x,y)=A_i(x,y)$ (for $i\in[1,s]$).
Thus $f_i\colon\Omega\to X\setminus\{i\}$.
Every fibre of every $f_i$ has size
$s-1$. Moreover, for distinct $i,k\in X$, and for
$j\ne i$, $h\ne k$,
\[
 \bigl|\{u\in\Omega:f_i(u)=j,\ f_k(u)=h\}\bigr|
 =
 \begin{cases}
  0, & (j,h)=(k,i)\text{ or }j=h,\\
  1, & \text{otherwise}.
 \end{cases}
\]
For $i,k\ne0$, this is precisely the orthogonality of $A_i$
and $A_k$. When one of $i,k$ is $0$, it follows from the row
conditions for an array-bic square.

Now define
\begin{align*}
 \mathcal P &= \{P_i,P_i^*:i\in X\} \cup \{P_u:u\in\Omega\}\\
 \mathcal L &= \{\ell_0,\ell_0^*\} \cup \{\ell_{ij}:i,j\in X,\ i\ne j\}.
\end{align*}
To streamline our notation, we declare $\ell_0:=\{P_i:i\in X\}$,
$\ell_0^*:=\{P_i^*:i\in X\}$, and 
\[
\ell_{ij}:=\{P_i,P_j^*\} \cup \{P_u:f_i(u)=j\}.
\]

Since every fibre of $f_i$ has size $s-1$, every line has
$s+1$ points. Also, $P_i$ lies on $\ell_0$ and the $s$
lines $\ell_{ij}$ with $j\ne i$, while $P_i^*$ lies on
$\ell_0^*$ and the $s$ lines $\ell_{ji}$ with $j\ne i$.
Finally, $P_u$ lies on the $s+1$ lines
$\ell_{i,f_i(u)}$, where $i\in X$.
Hence every point also lies on $s+1$ lines.

It remains to check intersections. The lines $\ell_0$ and
$\ell_0^*$ are disjoint, while $\ell_{ij}$ meets them in
$P_i$ and $P_j^*$, respectively. For two distinct lines
$\ell_{ij}$ and $\ell_{kh}$, we have
\[
 \ell_{ij}\cap\ell_{kh}
 =
 \begin{cases}
  \{P_i\}, & i=k,\\
  \{P_j^*\}, & j=h,\\
  \varnothing, & (k,h)=(j,i),\\
  \{P_u\}, & \text{otherwise},
 \end{cases}
\]
where in the last case $u$ is the unique element satisfying
$f_i(u)=j$ and $f_k(u)=h$.

Thus any two lines meet in at most one point, and every line has a
unique disjoint line, namely $\ell_0^*$ for $\ell_0$,
$\ell_{ji}$ for $\ell_{ij}$. Hence $(\mathcal P,\mathcal L)$ is a partial linear space. Since
$|\mathcal P|=2(s+1)+s(s-1) =s^2+s+2$, the $s+1$ lines through any point contain
$1+s(s+1)=s^2+s+1$ points. Therefore every point has a unique non-collinear point.
Consequently, $(\mathcal P,\mathcal L)$ is a regular antipodal plane of order $s$.

Finally, $P_i$ and $P_i^*$ are non-collinear, so they are
antipodal. Moreover, $P_{xy}=\ell_{0x}\cap\ell_{y0}$,
and the line through $P_i$ and $P_{xy}$ is
$\ell_{i,A_i(x,y)}$. Thus the construction in the forward
direction recovers $A_1,\ldots,A_s$, completing the
equivalence.
\end{proof}

By computational depth-first search, there is no $\MOAS(4,2)$.
Similarly, we could find a $\MOAS(5,2)$, but a $\MOAS(5,3)$ does not exist, 
and we could find a $\MOAS(6,2)$, but there is no $\MOAS(6,3)$.
There exists a $\MOAS(7,2)$, but a backtrack search with \textsf{Julia} (273 minutes on a MacBook)
shows that there is no $\MOAS(7,3)$.
This also shows there are no regular antipodal planes for these small orders, and that in some sense, we do not even get close to such a configuration. The \textsf{Julia} code is available upon request from the first author.

\begin{figure}[!ht]
$
\begin{bmatrix}
\cdot&2&5&3&4\\
0&\cdot&4&5&3\\
0&4&\cdot&2&5\\
0&5&3&\cdot&2\\
0&3&2&4&\cdot
\end{bmatrix},
\begin{bmatrix}
\cdot&0&4&5&3\\
1&\cdot&5&3&4\\
5&0&\cdot&4&1\\
3&0&1&\cdot&5\\
4&0&3&1&\cdot
\end{bmatrix}
$
\caption{A $\MOAS(5,2)$.}
\end{figure}

\begin{figure}[!ht]
$
\begin{bmatrix}
\cdot&2&5&6&3&4\\
0&\cdot&3&4&5&6\\
0&6&\cdot&2&4&5\\
0&5&2&\cdot&6&3\\
0&4&6&3&\cdot&2\\
0&3&4&5&2&\cdot
\end{bmatrix},
\begin{bmatrix}
\cdot&0&3&4&5&6\\
1&\cdot&4&3&6&5\\
6&0&\cdot&5&1&4\\
5&0&6&\cdot&3&1\\
4&0&1&6&\cdot&3\\
3&0&5&1&4&\cdot
\end{bmatrix}
$
\caption{A $\MOAS(6,2)$.}
\end{figure}

\begin{figure}[!ht]
$
\begin{bmatrix}
\cdot&6&4&5&2&3&7\\
0&\cdot&3&7&5&6&4\\
0&2&\cdot&4&7&5&6\\
0&3&7&\cdot&6&2&5\\
0&4&6&3&\cdot&7&2\\
0&7&5&2&4&\cdot&3\\
0&5&2&6&3&4&\cdot
\end{bmatrix},
\begin{bmatrix}
\cdot&0&5&4&6&7&3\\
7&\cdot&1&5&3&4&6\\
1&0&\cdot&7&4&6&5\\
3&0&6&\cdot&7&5&1\\
4&0&3&6&\cdot&1&7\\
5&0&7&3&1&\cdot&4\\
6&0&4&1&5&3&\cdot
\end{bmatrix}
$
\caption{A $\MOAS(7,2)$.}
\end{figure}

\section{Conclusion}

We have shown that there are interesting combinatorial objects associated to regular antipodal planes. The first and most well-known object is a biplane, and the existence of an example with more than 79 points is a long-standing open problem in design theory.
We have shown that a regular antipodal plane gives rise to an equidistant permutation array of index 1, and we do not yet have
good bounds on the size of such an object. Such a bound would resolve Payne's Conjecture. 
We also introduced sharply deranged sets and mutually orthogonal array-bic squares, that in some sense capture how improbable it is for a regular antipodal plane to exist, of order greater than 14. This leads us to the following questions arising from this work:
\begin{enumerate}[(Q1)]
\item Does there exist a regular antipodal plane of order greater than 14? In particular does there exist one of order 15?
\item Is it true that the maximum size $R(n,1)$ of an equidistant permutation array on $n$ points, and of index 1, is at most $\frac{1}{4}n^2+2$ for all $n\ge 3$?
\item Does there exist a $\MOAS(s,3)$ for $s>7$?
\end{enumerate}

\section*{Acknowledgements}

Sam Adriaensen is supported by grant number 12A3Y25N of Research Foundation - Flanders (FWO).
He is grateful to the University of Western Australia and John Bamberg in particular for their hospitality during his stay in Perth.
Jason Williford was supported in part by the Simons Foundation Collaboration Grant 711898. 
He thanks the Department of Mathematics at the University of Western Australia, and John Bamberg in particular, for their warm hospitality and for providing a stimulating research environment during his sabbatical visit.

The authors would like to thank Gordon Royle for sharing his adjacency matrices for the known biplanes. 

The computer algebra systems \textsf{Magma} \cite{Magma}, \textsf{GAP} \cite{GAP} and the programming language \textsf{Julia} \cite{Julia} were used for computations in this paper.

\end{document}